\documentclass[12pt,reqno]{amsart}

\usepackage{amssymb}
\usepackage{color}

\newtheorem{theorem}{Theorem}

\theoremstyle{remark}

\title[Monochromatic solution submitted to arXiv]
  {Monochromatic solutions to $x+y=z^2$ in the interval $[N,cN^4]$}

\author[P\'eter P\'al Pach]{P\'eter P\'al Pach$^\dag$}
\thanks{${}^\dag$ Partially supported by the National Research, Development and Innovation Office NKFIH
  (Grant Nr.~PD115978) and the J\'anos Bolyai Research
  Scholarship of the Hungarian Academy of Sciences. The author has also received funding from the European Research Council (ERC) under the European Union’s Horizon 2020 research and innovation programme (grant agreement No 648509). This publication reflects only its author's view; the European Research Council Executive Agency is not responsible for any use that may be made of the information it contains.}
\email{ppp@cs.bme.hu}
\address{Department of Computer Science and Information Theory, Budapest
  University of Technology and Economics, 1117 Budapest, Magyar tud\'osok
  k\"or\'utja 2., Hungary,  \newline 
  Department of Computer Science and DIMAP, University of Warwick, Coventry CV4 7AL, UK.}

\begin{document}

\baselineskip=16pt

\maketitle


\begin{abstract}
Green and Lindqvist  proved that for any 2-colouring of $\mathbb{N}$, there are in\-fi\-ni\-tely many monochromatic solutions to $x+y=z^2$. In fact, they showed the existence of a monochromatic solution in every interval $[N,cN^8]$ with large enough $N$. In this short note we give a different proof for their theorem and prove that a monochromatic solution exists in every interval $[N,10^4N^4]$ with large enough $N$. A 2-colouring of $[N,(1/27)N^4]$ avoiding monochromatic solutions to $x+y=z^2$ is also presented, which shows that in $10^4N^4$ only the constant factor can be reduced.
\end{abstract}

\section{Introduction}

Csikv\'ari, Gyarmati and S\'ark\"ozy \cite{csgys} proved that the equation $x+y=z^2$ is not partition regular. Indeed, they gave a 16-colouring of $\mathbb{N}$ with no monochromatic solutions
to $x + y = z^2$ other than the trivial one $x = y = z = 2$. Recently, Green and Lindqvist \cite{gl} have shown that such a colouring also exists with only 3 colours, but for any 2-colouring of $\mathbb{N}$ there are infinitely many monochromatic solutions to $x+y=z^2$. In fact, from their proof it also follows that for every large enough $N$ there is a monochromatic solution in the interval $[N,cN^8]$ (where $c$ is a huge constant that can be explicitly given). The proof of Green and Lindqvist uses several tools from additive combinatorics.

In this paper we give another, shorter proof for this theorem. The proof is elementary and involves several combinatorial ideas. Our result is the following:

\begin{theorem}\label{thm}
Every 2-colouring of $\mathbb{N}$ has infinitely many monochromatic solutions to $x+y=z^2$. Moreover, for every large enough $N$ there is a monochromatic solution in $[N,10^4N^4]$.
\end{theorem}

This strengthens the result of Green and Lindqvist, in the sense that it verifies the existence of a monochromatic solution in a much shorter interval, $[N,10^4N^4]$ instead of $[N,cN^8$]. Note that the exponent 4 in $10^4N^4$ can not be further improved, since colouring  $[N,N^2/3]$ with the first colour and $(N^2/3,N^4/27]$ with the second colour avoids any monochromatic solution for $x+y=z^2$.

The proof of Theorem~\ref{thm} is given in Section 2. Throughout the paper we use the notion $[n]:=\{1,2,\dots,n\}$ and by an interval $[a,b]$ we mean the set of integers that are at least $a$ and at most $b$.

Finally, we shall mention two related works. Khalfallah and Szemer\'edi \cite{ksz} have shown that for any finite colouring of $\mathbb{N}$ there is a monochromatic solution
to $x+ y = z^2$ with $x$ and $y$ having the same colour (but not necessarily $z$).

Lindqvist \cite{l} considered the modular version and showed that
if $p > p_0(k)$ is a prime and if $\mathbb{Z}/p\mathbb{Z}$ is $k$-coloured, then there are $\gg_k p^2$ monochromatic solutions to $x + y = z^2$.

\section{Proof}

It suffices to prove the second statement of Theorem~\ref{thm}.

Let $c:[N,10^4N^4]\to \{-1,1\}$ be an arbitrary 2-colouring.

If all elements of $[9N,80N^2]$ are coloured with the same colour, then $x=N^2,y=80N^2,z=9N$ is a monochromatic solution.

Otherwise, let $k\in [9N,80N^2]$ be such that $c(k)\ne c(k+1)$. Without loss of generality we shall assume that $c(k)=1$ and $c(k+1)=-1$.

If there exists some $i\in [N,k^2-N]$ with $c(i)=c(k^2-i)=1$, then $x=i,y=k^2-i,z=k$ is a monochromatic solution. Therefore, we can assume that 
\begin{equation}\label{eq1}
c(i)+c(k^2-i)\leq 0
\end{equation}
holds for every $i\in [N,k^2-N]$.

Similarly, if there exists some $i\in [N,(k+1)^2-N]$ with $c(i)=c((k+1)^2-i)=-1$, then $x=i,y=(k+1)^2-i,z=k+1$ is a monochromatic solution. Therefore, we can also assume that 
\begin{equation}\label{eq2}
c(i)+c((k+1)^2-i)\geq 0
\end{equation}
holds for every $i\in [N,(k+1)^2-N]$.

Now, let $j\in [N,k^2-N]$. By taking $i=j$ in \eqref{eq1} and $i=j+2k+1$ in \eqref{eq2} we obtain
\begin{equation}\label{eq3}
c(j)+c(k^2-j)\leq 0
\end{equation}
and
\begin{equation}\label{eq4}
c(j+2k+1)+c(k^2-j)\geq 0.
\end{equation}
Inequalities \eqref{eq3} and \eqref{eq4} yield that $c(j)\leq c(j+2k+1)$ holds for every $j\in [N,k^2-N]$. That is, for every $j\in [N,N+2k]$ we have
\begin{equation}\label{monoton}
c(j)\leq c(j+(2k+1))\leq c(j+2(2k+1))\leq\dots\leq c(j+m_j(2k+1)),
\end{equation}
where $m_j$ is the largest integer such that $j+(m_j-1)(2k+1)\leq k^2-N$.
Note that $j+m_j(2k+1)> k^2-N$.
For $j\in [N,N+2k]$ let $H_j:=\{j,j+(2k+1),\dots,j+m_j(2k+1)\}$.

We obtained that restricting the colouring to any mod $2k+1$ residue class $H_j$ it is monotonic, in the sense that the pattern of colours looks like $-1,-1,\dots,-1,1,1,\dots,1$. Let us introduce a function which tells us where the breaking point is.

For $j\in [N,N+2k]$ let $f(j)=\infty$, if all elements of $H_j$ are coloured $-1$, otherwise let $f(j)$ be the smallest element of $H_j$ which is coloured $1$. That is, $f(j)=j+l(2k+1)$, if $c(j)=\dots=c(j+(l-1)(2k+1))=-1$ and $c(j+l(2k+1))=\dots=c(j+m_j(2k+1))=1$. (If $c(j)=\dots=c(j+m_j(2k+1))=-1$, then $f(j)=\infty$.)

The function $f$ is defined on a complete residue system modulo $2k+1$. Let us extend it to the set of integers, in such a way that for an integer $j'$ let $f(j')=f(j)$ with the unique $j\in[N,N+2k]$ which is congruent with $j'$ modulo $2k+1$. Similarly, $H_{j'}:=H_j$ for the unique $j\in [N,N+2k]$ congruent with $j'$ modulo $2k+1$. 

Assume now that for some $j$ we have $f(j)+f(k^2-j)\leq k^2$.
By taking the sum of an element from $H_j$ and an element from $H_{k^2-j}$ we always get a number which is congruent with $k^2$ modulo $2k+1$. Note that those numbers in $H_j$ which are coloured 1 are next to each other, and so do those numbers in $H_{k^2-j}$ that are coloured 1. Therefore, the set of such possible sums is also an interval (in the mod $2k+1$ residue class of $k^2$). The smallest such sum is $f(j)+f(k^2-j)$ and the largest one is at least $2(k^2-N)$. As $k^2\in [f(j)+f(k^2-j),2(k^2-N)]$, from the classes  $H_j$ and $H_{k^2-j}$ we can choose two elements coloured $1$ whose sum is $k^2$, and we obtain a monochromatic solution, since $c(k)=1$.

From now on, let us assume that $f(j)+f(k^2-j)> k^2$ for every $j$. Note that this implies that for every $j$ we have $f(j)+f(k^2-j)\geq k^2+2k+1$. Let $A=\{j\in [N,N+2k]: f(j)\geq (k^2+2k+1)/2\}$. Since $f(j)+f(k^2-j)\geq k^2+2k+1$ for every $j$, we have $|A|\geq k+1$. By the pigeon-hole principle $A+A$ contains elements from all of the residue classes modulo $2k+1$. 



Let $m\in[0.2k,0.8k]$ be an integer. (Note that $m\geq 0.2k\geq N$.) We can choose $j_1,j_2\in A$ such that $m^2\equiv j_1+j_2 \pmod{2k+1}$.

For $j\in A$ let $g(j)$ denote the largest element from $H_j$ which is coloured $-1$. That is, if $f(j)<\infty$, then $g(j)=f(j)-(2k+1)$ and if $f(j)=\infty$, then $g(j)=\max(H_j)$. Note that according to $j\in A$, we have $g(j)\geq (k^2-2k-1)/2$, since either $g(j)=f(j)-(2k+1)\geq (k^2-2k-1)/2$ or $g(j)=\max(H_j)>k^2-N>k^2-k$. 

From the residue class of $m^2$ we can write each element between $j_1+j_2$ and $g(j_1)+g(j_2)$ as a sum of two numbers coloured $-1$ (taken from $H_{j_1}$ and $H_{j_2}$). Since $j_1+j_2\leq 2\cdot(N+2k)<6k<(0.2k)^2$ and  $g(j_1)+g(j_2)\geq k^2-2k-1>(0.8k)^2$, the number $m^2$ is also the sum of two integers which are coloured $-1$. If $c(m)=-1$, then this results a monochromatic solution with $z=m$. 
Therefore, we can assume that all integers from $[0.2k,0.8k]$ are coloured 1, and so do all the elements from their mod $2k+1$ residue classes up to $k^2-N>k^2-k$, according to \eqref{monoton}. Note that $k^2\equiv 1.5k+1 \text{ or } 0.5k+0.5 \pmod{2k+1}$ (depending on the parity of $k$). Since in both cases the modulo $2k+1$ residue of $k^2$ can  be expressed as a sum of two residues taken from $[0.2k,0.8k]$ we obtain a monochromatic solution with $z=k$.

{\bf Acknowledgements.} The author would like to thank P\'eter Csikv\'ari and Hong Liu for reading an earlier version of this note and for their useful comments.


\vfill
\begin{thebibliography}{McWS77}


\bibitem[1]{csgys}
P.~Csikv\'ari, A.~S\'ark\"ozy, K.~Gyarmati:
\newblock {\em Density and Ramsey type results
on algebraic equations with restricted solution sets
},
\newblock Combinatorica 32 (2012), 425--449.



\bibitem[2]{gl}
B.~Green, S.~Lindqvist:
\newblock {\em Monochromatic solutions to  $x+y=z^2$},
\newblock Canadian Journal of Mathematics, to appear

\bibitem[3]{ksz}
A.~Khalafallah, E.~Szemer\'edi:
\newblock {\em On the Number of Monochromatic Solutions
of $x+y=z^2$},
\newblock Combinatorics, Probability and Computing 15 (2006), (15) 1-2,
213--227.

\bibitem[4]{l}
S.~Lindqvist:
\newblock {\em Monochromatic solutions to $x + y$ a square in $\mathbb{Z}/q\mathbb{Z}$},
\newblock available at
http://people.maths.ox.ac.uk/lindqvist/notes/xysumsquare.pdf



\end{thebibliography}
\end{document}